\documentclass[12pt]{article}

\usepackage{amsmath, amsthm, amsfonts, makeidx}
\DeclareMathSymbol{\subsetneq}{\mathord}{AMSb}{"26}

\newtheorem{lemma}{Lemma}[section]

\newtheorem{theorem}[lemma]{Theorem}

\newtheorem{proposition}[lemma]{Proposition}
\newtheorem{corollary}[lemma]{Corollary}
\theoremstyle{definition}

\newtheorem{definition}[lemma]{Definition}
\newtheorem{example}[lemma]{Example}

\newtheorem{remark}[lemma]{Remark}

\newtheorem{assumption}[lemma]{Assumption}

\newcommand{\lp}{\longrightarrow}

\newcommand{\mb}{\mathbb}

\newcommand{\G}{\mathcal{G}}
\newcommand{\F}{\mathcal{F}}

\newcommand{\C}{\mb{C}}

\newcommand{\Z}{\mb{Z}}
\newcommand{\N}{\mb{N}}
\newcommand{\Q}{\mb{Q}}

\newcommand{\DER}{\textup{DER}}

\newcommand{\grad}{\operatorname{grad}}
\renewcommand{\ker}{\operatorname{ker}}

\renewcommand{\deg}{\operatorname{deg}}
\newcommand{\trdeg}{\textup{trdeg}}
\newcommand{\LND}{\textup{LND}}

\newcommand{\gr}{\textup{gr}}
\newcommand{\GR}{\textup{GR}}

\newcommand{\kk}{\mathbb{C}}
\newcommand{\KK}{\mathbb{C}}
\newcommand{\ff}{\tilde{f}}
\newcommand{\hh}{\tilde{h}}
\newcommand{\NN}{\mathcal{N}}
\newcommand{\ml}{\textup{ML}}

\title{Rigid rings and Makar-Limanov techniques}
\author{Anthony J. Crachiola\\
\small Department of Mathematical Sciences\\
\small Saginaw Valley State University\\
\small 7400 Bay Road\\
\small University Center, Michigan 48710, USA\\
\small crachiola@member.ams.org\\ \ \\
Stefan Maubach\footnote{Funded by Veni-grant of council for the
physical sciences, Netherlands Organisation for scientific research
(NWO).}\\ \small Radboud University Nijmegen\\\small Toernooiveld 1,
The Netherlands\\ \small s.maubach@math.ru.nl}

\begin{document}

\maketitle

\begin{abstract}
A ring is rigid if there is no nonzero locally nilpotent derivation
on it. In terms of algebraic geometry, a rigid coordinate ring
corresponds to an algebraic affine variety which does not allow any
nontrivial algebraic $\G_a$-action. Even though it is thought that
``generic'' rings are rigid, it is far from trivial to show that a
ring is rigid. In this paper we provide several examples of
rigid rings and we outline
two general strategies to help determine if a ring is rigid, which
we call
``parametrization techniques'' and
``filtration
techniques''.
We provide many little tools and lemmas which may be useful in other
situations. Also, we point out some pitfalls to beware when using
these techniques. Finally, we give some reasonably simple
hypersurfaces for which
the question of rigidity remains unsettled.
\end{abstract}

\section{Introduction}




%
%



In our article, all rings are commutative and contain 1. We will use
$\N$ to denote the set of natural numbers, $\Q$ to denote the field
of rational numbers, and $\Q^+$ to denote the set of positive
rational numbers.

If $A$ is a ring, a derivation $D:A\lp A$ is a linear map satisfying
the Leibniz rule $D(ab)=aD(b)+D(a)b$ for all $a,b\in A$. This
derivation is called locally nilpotent if for each $a\in A$ there
exists $n\in \N$ such that $D^n(a)=0$. We write $\DER(A)$ for the
set of derivations on $A$, and $\LND(A)$ for the set of locally
nilpotent derivations on $A$. We use LND to abbreviate locally nilpotent derivation. Let $A^D := \ker(D)=\{a\in A ~|~
D(a)=0\}$. We define a ring to be a {\bf rigid ring} if it admits no
LND except the zero derivation. Similarly, an affine variety $X$
is called rigid if $\mathcal{O}(X)$ is rigid.

We will focus on the field $\kk$, though many things will apply to any field of characteristic zero. 
We denote $\kk^{[n]}:=\kk[X_1,\ldots,X_n]$ for the polynomial
ring in $n$ variables.
If $X$ is a variable of $\kk^{[n]}$, then
$\partial_X:=\frac{\partial}{\partial X}$ is an LND. In this paper we focus on $\kk$-algebras, and by a
``derivation'' on a $\kk$-algebra, we implicitly mean a derivation
which is $\kk$-linear.

Our approach will be primarily algebraic. However, there is an
important geometric interpretation of the subject matter.
%
Any LND $D$ of the $\kk$-algebra $A$ gives rise to an algebraic
action of the additive group $\G_{a}$ of $\kk$ on $A$ via
\begin{equation*}
\exp (tD)(a) = \sum
\limits_{i=0}^{\infty}\frac{t^{i}}{i!}D^{i}(a)
\end{equation*}%
for $t\in \kk$ and all $a\in A.$ Conversely, an algebraic action
$\sigma$ of
$\G_{a}$ on $A$ yields an LND on $A$ via%
\begin{equation*}
\frac{\sigma (t,a)-a}{t}\,\Bigr|_{t=0}.
\end{equation*}%
The kernel  $A^{D}$ coincides with the ring of
invariants of the corresponding $\G_{a}$-action.

The
Makar-Limanov invariant of $A$, denoted $\ml(A)$, is defined as the
intersection of the kernels of all LNDs on
$A$. Thus to show that $A$ is rigid means to show that $\ml(A)=A$.
The Makar-Limanov invariant has been used with success because there
are techniques to compute it (also due to Makar-Limanov). In some
sense, these techniques
are more valuable than
%
the invariant itself. (Note however, that these techniques are not
omnipotent -- there's no algorithm to compute the ML invariant, for
one.)
In this paper we will focus these techniques onto the task of
showing that certain rings are rigid.




\section{Derivations on non-domains}

There are advantages to working with a domain (or an irreducible algebraic variety) when studying LNDs and rigidity. Often, however, we are lead to consider LNDs on non-domains.

\begin{proposition}\label{ND-one}
Suppose $\G_a$ acts algebraically on a reducible affine variety $X$.
Then the $\G_a$-action restricts to each irreducible component of
$X$, and to each of their intersections.
\end{proposition}

\begin{proof}
Let $X=X_1\cup X_2\cup \cdots \cup X_r$ be the decomposition of $X$
into irreducible components. The group $\G_a$ acts algebraically on
$X$. Let $\varphi \in \G_a$, which is an automorphism of $X$. It
thus sends irreducible factors to irreducible factors, and thus
permutes the set $\{X_1,\ldots,X_r\}$. So for each $\varphi\in
\G_a$ we get a permutation $\sigma(\varphi)\in \textup{Sym}(r)$.
Thus we have an induced algebraic map $\G_a\overset{\sigma}{\lp}
\textup{Sym}(r)$. Now $\ker(\sigma)$ is a subvariety of $\G_a\cong
\C$, hence it is either a finite subset of $\C$ (which is impossible) or
equal to $\C$. So $\ker(\sigma)=\G_a$, which means that
$\G_a(X_i)=X_i$ for each $1\leq i\leq r$.
\end{proof}

\begin{corollary} \label{ND-C}
Let $f\in \C^{[n]}$ be a radical polynomial, i.e. $f=f_1f_2\cdots f_r$ where each $f_i$
 is irreducible and different from the others.
Write $A:=\C^{[n]}/(f)$ and let $\bar{f}_i$ be the image of $f_i$ in $A$.
If $D\in \LND(A)$ then $D(\bar{f}_i)\subseteq \bar{f}_iA$ for all $1\leq i \leq r$.
\end{corollary}

This follows immediately from Proposition \ref{ND-one}, as apparently the
$\G_a$-action induced by $D$ restricts to the domain
$\C^{[n]}/(f_i)$, meaning that $D(\bar{f}_i)\subseteq \bar{f}_iA$.

\begin{corollary}\label{NC-C2}
%
Let $X$ be an affine variety with irreducible components $X_1, \ldots, X_n$. If each $X_i$ is rigid, then $X$ is rigid. Furthermore, there is a 1-1 correspondence between the algebraic $\G_a$-actions on $X$ and  sets of algebraic actions $\G_a\times X_i \overset
{\varphi_i}{\lp} X_i$ such that $\varphi_i\arrowvert_{X_i\cap X_j}=
\varphi_j \arrowvert_{X_i\cap X_j}$ for each $i,j$.

\end{corollary}

An application of the above:

%

\begin{lemma}\label{XYrigid}
Let $R$ be a finitely generated domain over $\C$. Then $R$ is rigid if and only if $R[X,Y]/(XY) = R[x,y]$ is rigid.\end{lemma}

\begin{proof}
If $R$ is not rigid then any nonzero LND of $R$ extends trivially to a nonzero LND of $R[x,y]$.
Suppose now that $R$ is rigid and let $D$ be an LND of $R[x,y]$. Then $D(x) \subseteq xR[x,y]$ and
$D$ restricts to an LND $D_x$ of $R[x,y]/(y) \cong R[x]$, where $x$ is a free variable over $R$.
By Theorem 3.1 in \cite{tony-ml}, we have $\ker(D_x) = R$. Thus
$D_x = a \partial_x$
for some $a \in R$. 
Now $a =D_x (x) \in xR[x,y]/(y)$ which means
%
that $a \in R \cap xR[x] = \{ 0 \}$ and thus $D_x = 0$. Likewise, the restriction of $D$ to $R[x,y]/(x)$ is zero. Hence $D = 0$.
\end{proof}

\section{Parametrization techniques}


Let $A$ be a domain containing $\mathbb{Q}$. Let $D \in \LND(A)$
with $D \ne 0$. Let us summarize some well known and commonly used
facts (cf. Chapter 1 of \cite{Essenboek}). There exists $p \in A$
such that $ D^{2}(p)=0$ while $D(p)\not=0$. Set $q:=D(p)$ and $s  =
p/q$. There is a unique extension of $D$ to a locally nilpotent
derivation (which we still denote by $D$) of $\tilde{A}:=A[q^{-1}]$
with $\tilde{A}^{D} = A^D[q^{-1}]$ and $D(s) = 1$. Moreover,
$\tilde{A} = A^D[q^{-1}][s]$ and $D$ is the usual derivative
$\partial_{s} = \partial/\partial s$ on $\tilde{A}$.  Now $D$
extends uniquely to $k[s]$, where $k$ is the quotient field of
$\tilde{A}^D$. In fact, $D$ also extends uniquely to $K[s]$, where
$K$ is the algebraic closure of $k$. This brings us to the starting
point for parametrization techniques.

\begin{lemma} (Parametrization Lemma)\label{PT1} 
\label{L1} Let $A$ be a domain containing $\mathbb{Q}$. Let $D \in
\LND(A)$ with $D \ne 0$. Then $A$\ embeds into a polynomial ring
$K[S]$, where $K$ is some algebraically closed field,
and $D$ is the restriction of $\partial_{S} \in \LND(K[S])$ to $A$.
\end{lemma}

The Parametrization Lemma
allows us to reduce questions on LNDs to
questions on polynomials, and to use obvious facts about polynomials
to discover less obvious facts about LNDs.
The following lemma contains well known facts (see for instance Proposition
1.3.32 of \cite{Essenboek} or Section 1.4 of \cite{Fre06}). We
include the proof to illustrate how parametrization can be applied.

%
%

%

\begin{lemma} \label{BS1} 
Let $A$ be a domain containing $\Q$. Let $D\in \LND(A)$. \\
{\bf a)} If $u \in A$ is a unit, then $D(u)=0$.\\
{\bf b)} $A^D$ is factorially closed, that is: if $ab\in A^D, a,b\in A$,  then $a,b,\in A^D$.\\
{\bf c)} $A^D$ is algebraically closed in $A$.\\
{\bf d)} If $\lambda,a\in A$, $a\not =0$, such that $D^n(a)=\lambda a$ for some $n\geq1$, then $\lambda=0$.\\
{\bf e)} If $a \in A, a\not =0$ and $E\in \DER(A)$ such that $D=aE$, then $E\in \LND(A)$ and $a\in A^D=A^E$.\\
{\bf f)} If, moreover, $A$ is a $\kk$-algebra with $\trdeg_{\kk}(A)
< \infty$ and $D \not =0$, then
$\trdeg_{\kk}(A^D)=\trdeg_{\kk}(A)-1$.
\end{lemma}

\begin{proof}
Suppose $D(ab)=0$. Write $a,b\in K[S]$ using notations of Lemma
$\ref{PT1}$. Then $\deg_S(ab)=0$, hence both $a$ and $b$ must be
constant, and thus $D(a)=\partial_S(a)=0$ and
$D(b)=\partial_S(b)=0$. This proves b). Since $D(1)=0$, we get a) as
a consequence of b). Part c) also follows immediately from b). Now
assume $D^n(a)=\lambda a$ and again write $a\in K[S], D=\partial_S$.
Since $\deg_S(\partial^n_S(a)) \leq \deg_S(a) - n$,
it is clear that this can only be true if $n=0$ or if $\lambda =0$,
proving d). To prove e), consider $D(a)=aE(a)$. By d), $E(a)=0$. For
any $b\in A$ we find $n\in \N$ such that $0=D^n(b)=(aE)^n(b)$. Since
$E(a)=0$, $a$ commutes with $E$, i.e. $E(ac)=aE(c) $ for all $c\in
A$. Thus $0=a^nE^n(b)$, and since $a\not =0$ we get $E^n(b)=0$. Thus
$E$ is locally nilpotent. Part f) follows from the earlier
observation that $A[q^{-1}] = A^D[q^{-1}][s]$, where $s = p/q$ for
some $p \in A$ and $q = D(p) \in A^D$.
\end{proof}


\begin{lemma}\label{magic1} ($\Q^+$-twist commutator lemma)
Let $f,g\in \kk[X]$, both nonzero. Suppose that $ f'g=h fg'$ where
$h \in \kk[X]$ and $h \not = 0$. Then  either $f,g\in \kk$
simultaneously or $h\in \Q^+$.
\end{lemma}

\begin{proof}
It is clear that $f$ and $g$ can only be a constant simultaneously.
So assume $\deg(f)=n>0, \deg(g)=m>0$. Since $\deg(f'g) = \deg(fg') =
n+m-1 > 0$, we find that $h$ is a constant. Then comparing the
highest degree terms of $f'g$ and $fg'$ we get that $n=h m$, hence
the result follows.
\end{proof}


\begin{corollary}\label{magic} ($\Q^+$-twist commutator lemma)
Let $A$ be a domain containing $\Q$. Let $D\in \LND(A)$. Let $a, b,
\lambda \in A$ and $\mu\in A^D$, all nonzero. If $ \lambda aD(b)=\mu
bD(a)$ then either $D(a)=D(b)=0$ or $\lambda\mu^{-1}\in \Q^+$.
\end{corollary}

\begin{proof}
Follows immediately from Lemmas \ref{PT1} and \ref{magic1}.
\end{proof}

\begin{lemma}\label{Lefschetz}
Let $K \supset \kk$ be a field. Let $P \in \kk^{[n]}$. If
$P(f_1,\ldots,f_n) \ne 0$ for all $(f_1,\ldots, f_n) \in \kk[S]^n
\setminus \kk^n$, then $P(f_1,\ldots,f_n) \ne 0$ for all
$(f_1,\ldots, f_n) \in K[S]^n \setminus K^n$.
\end{lemma}


\begin{proof}
An application of the Lefschetz Principal.
\end{proof}



\begin{corollary} \label{PT2}
Let $A=\kk^{[n]}/(P)$ where $P\in \kk^{[n]}$ is irreducible.
Suppose that 
$P(f_1,\ldots,f_n) \ne 0$ for all $(f_1,\ldots, f_n) \in \kk[S]^n
\setminus \kk^n$.
Then $A$ is rigid.
\end{corollary}

%

\begin{proof}
Write $A =  \kk[x_1, \ldots, x_n]$, where $x_i = X_i + (P)$.
Suppose
$D\in \LND(A)$ is nonzero. By Lemma \ref{PT1} we view each $x_i$ as
an element of $K[S]$. By assumption, each $x_i$ must be constant as
a polynomial in $S$. But that means, that the extension
$\frac{\partial}{\partial S}$ of $D$ on $K[S]$ is zero on each
$x_i$, i.e. $\frac{\partial}{\partial S}(x_i)=D(x_i)=0$ for all $i$.
But then $D=0$, contradiction.
Hence there exists no nonzero $D\in \LND(A)$, so $A$ is rigid.
\end{proof}

We now summarize the lemmas that are typically useful for
parametrization situations. The nicknames are chosen in reference to
Mason's Theorem \ref{Mason}, discussed below.
Mini-Mason can also be found in a certain form as Lemma 9.2 in
\cite{Fre06}, and Lemma 2 in \cite{lenny}. These quoted lemmas are
excellent examples of what is essentially the parametrization
technique.

\begin{lemma}\label{MiniMason} \label{PT3}
 Let $f,g\in K[S]$ where $K$ is an algebraically
closed field of characteristic zero. \\
{\bf a)}
(`` Mini-Mason'') If $f^a+g^b \in K^*$ where $a,b\geq 2$, then $f,g\in K$.\\
{\bf b)}
(``Extended Mini-mason'') If  $f^a+g^bQ(g)\in K^*$ where
$Q(T)\in K[T]$ is nonzero, and where $a,b,\in \N$ such that
$\deg_T(Q)+1\leq (a-1)(b-1)$, then $f,g\in K$.\\
{\bf c)}
(``Twisted Mason'') If $f^ag^b+h^c\in K^*$ where $a,b,c\geq 2$,
then $f,g,h\in K$.
\end{lemma}

\begin{proof}
Part a) follows from b). For b),
note
that $\gcd(f,g)=1$. Differentiate  with
respect to $S$ to obtain
$af'f^{a-1}=-g^{b-1}(bg'Q(g)+gg'Q'(g))$.
Write $F$ (resp. $G$) for the degree of $f$ (resp. $g$) in $S$. Note
that
\[
\deg(bg'Q(g)+gg'Q'(g)) \leq (\deg_T(Q) + 1)G - 1.
\]
Since $\gcd(f,g)=1$, $f^{a-1}$ divides $bg'Q(g)+gg'Q'(g)$ and
$g^{b-1}$ divides $f'$. This means that $(a-1)F \leq
(\deg_T(Q)+1)G-1$ and $(b-1)G=F-1$. Combining these we obtain
$(a-1)(b-1)<\deg_T(Q)+1$ which contradicts the assumption exactly.

The proof of c) is similar to that of b). We have
$\gcd(f,h)=\gcd(g,h)=1$, and we take the derivative and get
$f^{a-1}g^{b-1}(af'g+bfg')=-ch'h^{c-1}$. From this we conclude that
$af'g + bfg'$ divides $h^{c-1}$ and $f^{a-1}g^{b-1}$ divides $h'$.
Write $F,G,H$ for the degree of $f,g,h$, respectively. We have
\[ F+G-1=(c-1)H, \quad (a-1)F+(b-1)G=H-1\]
and from this we obtain
\[ [(a-1)(c-1)-1]H+(b-a)G+a=0, \quad [(b-1)(c-1)-1]H+(a-b)F+b=0. \]
If $b\geq a$ the first equation is impossible, while if $a\leq b$
the second equation is impossible.
\end{proof}

Mason's Theorem, below, provides a very useful technique for
constructing rigid rings.
We will use the following generalization due to M. de Bondt. The
original theorem handles the case $n=3$. While R.~Mason did publish
a proof \cite{Mason}, the original proof is due to W.~Stothers
\cite{Stothers}. In the statement of the theorem, $\NN(h)$ denotes
the number of distinct zeroes of a polynomial $h\in K[S]$, where $K$
is algebraically closed.


\begin{theorem}[Theorems B1.3 and B1.4 of \cite{MichielMason}]
\label{Mason} Let $f_{1},f_{2},\ldots ,f_{n}\in K[S]$ where $K$ is
an algebraically closed field of characteristic 0. Assume not all
$f_i$ are constant, and
$f_{1}+f_{2}+\cdots +f_{n}=0$.
Additionally, assume that for every $1\leq i_{1}<i_{2}<\ldots <i_{s}\leq n$,
$f_{i_{1}}+f_{i_{2}}+\cdots +f_{i_{s}}=0\Longrightarrow
\gcd(f_{i_{1}},f_{i_{2}},\ldots ,f_{i_{s}})=1$.
Then both of the following inequalities hold.
%
\begin{align*}
\max\{\deg(f_1),\ldots,\deg(f_n)\} &< \dfrac{(n-1)(n-2)}{2}\NN(f_1f_2\cdots f_n)\\
\max\{\deg(f_1),\ldots,\deg(f_n)\} &< (n-2)(\NN(f_1)+\NN(f_2)+\cdots +\NN(f_n))
\end{align*}

\end{theorem}




%

\section{Examples using parametrization techniques}

In this section we consider examples of the form
$A:=\KK[X_1,\ldots,X_n]/(f)$. Also, we denote $x_i:=X_i+(f)$, i.e.
$A=\KK[x_1,\ldots,x_n]$.

\subsection{Rigidity of $X^aY^b-Z^c$}

Write $A:=\kk[X,Y,Z]/(X^aY^b-Z^c)$.
Trivial cases for which $A$ is not rigid occur when at least
%
one of $a,b,c$ equals zero
as well as when at least
one of $a, b, c$ equals 1. We will show that all other cases are rigid.

\begin{lemma} \label{xayb-zcdomain} Let $a,b,c\geq 2$ and $\gcd(a,b,c)=1$. Then $A$
is rigid. 
\end{lemma}

The proof follows directly from Twisted Mason \ref{MiniMason}
and the parametrization technique. We proceed to the general case:

\begin{theorem}\label{case1}
If $a,b,c\geq 2$ then $A$ is rigid.
\end{theorem}

\begin{proof}
Note that for $A$ to be
a domain means that $\gcd(a,b,c)=1$.
Lemma \ref{xayb-zcdomain} covers the case when $A$ is a domain. 
Let $d=\gcd(a,b,c)>1$, and denote $a'=a/d, b'=b/d, c'=c/d$. Now 
using Lemma~\ref{ND-C} we see that any $D \in \LND(A)$ restricts to the irreducible components, each of which is isomorphic to $X^{a'}Y^{b'}-Z^{c'}$. Denote the restriction to 
$A':=\C[x,y,z]:=\C[X,Y,Z]/(X^{a'}Y^{b'}-Z^{c'})$
by $\delta$.
The restriction of $\delta$ 
to the intersection of these irreducible components, which corresponds to the ring $\C[X,Y,Z]/(XY,Z)$, is zero by Lemma \ref{XYrigid}. We will show that the only LND on $A'$ which is zero modulo $(xy,z)$ is the zero derivation.
We distinguish several cases, where we state $a'\leq b'$.\\
Case 1. Suppose  $a',b', c'>1$. 
Then $A$ is rigid by
Lemma \ref{xayb-zcdomain}.\\
Case 2. Suppose  $a'=1$ and $b',c'>1$. By \cite{lml} any nonzero LND $\delta$ on $A'
$ has the same kernel, namely $\C[y]$, and is of the form 
$\delta = p(y)(c'z^{c'-1}\partial_x+y^{b'}\partial_z)$. It must be zero modulo $(xy,z)$, so $\delta(z)=p(y)y^{b'}\in (xy,z)$ implies that $p=0$. Thus $D$ is zero on each irreducible component, thus $D$ is zero itself.\\
Case 3. Suppose  $c'=1$. $D$ induces an LND $\delta$ on
$A'
\cong \C[X,Y]$
which is zero modulo the ideal $(XY)$. So in particular, $\delta(X)=\alpha XY$ and $\delta(Y)=\beta XY$ for some $\alpha,\beta\in \C[X,Y]$ which means by Lemma \ref{BS1} part d) that $\delta(X)=\delta(Y)=0$. \\
Case 4. suppose $a'=b'=1$ and  $c'>1$. The hardest case. By Proposition 2.8 of \cite{Daigle}
we may assume that $\ker(\delta)=\C[x^{-1}(z-\lambda x^q)^{c'}]$ for some $\lambda\in \C$. 
Define the automorphism $\varphi(x,y,z)=(x,x^{-1}(z-\lambda x^q)^{c'}, z-\lambda x^q)$
with inverse $\varphi^{-1}(x,y,z)=(x,x^{-1}(z+\lambda x^q)^{c'}, z+\lambda x^q)$.
Then $\Delta:=\varphi^{-1}\delta\varphi $ has as kernel $\C[y]$, and thus 
$\Delta=p(y)(c'z^{c'-1}\partial_x+y\partial_z)$. 
Now $\delta(x)=\varphi \Delta \varphi^{-1}(x)=\varphi(p(y)c'z^{c'-1})=p(\varphi(y))c'\varphi(z)^{c'-1}$. We assumed $\delta=0 \mod (xy,z)$ so $0=\delta(x)\mod(xy,z)$ equals
$p(\varphi(y))(-\lambda x^q)^{c'-1}$. This implies either $p=0$ (meaning that $\delta=0$) or $\lambda=0$. In the latter case $\delta=\Delta$, and $\Delta(z)=p(y)y\in (xy,z)$ only if $p=0$.  
 \end{proof}

\subsection{Rigidity of $X^dY+Z^dP(Y)$ with $\deg(P)\leq (d-1)^2$}

\begin{theorem} \label{EX2t}
Let $A:=\KK[X,Y,Z]/(X^dY+Z^dP(Y))=\kk[x,y,z]$, where $d\geq 2$,  
and
$P(T)\in \KK[T]$ with $P(0)\not=0$ and
$\deg(P)\leq (d-1)^2$. Then $A$ is rigid.
\end{theorem}


By Extended Mini-Mason \ref{MiniMason}, the degree restriction on
$Q$ implies that there are no nonconstant parametrizations of the
form $F^d + H^dQ(H) = 1$, where $F, H \in \kk^{[1]}$. The conclusion
then follows from Proposition \ref{EX2} below.
%
The below proposition does include more cases than the above
theorem, though (for example, if $Q=T^e$ where $e\geq 2$).
%
%
%

\begin{proposition} \label{EX2}
Let $P(T)\in \KK[T]$, $A:=\KK[X,Y,Z]/(X^dY+Z^dP(Y))=\kk[x,y,z]$
where we have $d\geq 2$, $P(0)\not=0$. Write $P(T)=TQ(T)+P(0)$.
Assume there exist no nonconstant parametrizations of the
form
$F^d + H^dQ(H) = 1$, where $F, H \in \kk^{[1]}$. Then $A$ is rigid.
\end{proposition}

\begin{proof}
Suppose $D$ is a nonzero LND on $A$. By the
Parametrization Lemma \ref{PT1} we view $A \subseteq K[S]$ for some
algebraically closed field $K$ with $D$ as the restriction of
$\frac{\partial}{\partial S}$. Write $x = f(S)$, $y = g(S)$, and $z
= h(S)$. Let $\alpha = \gcd(f,h)$ in $K[S]$ and write $f = \alpha
\ff$, $h = \alpha \hh$. From the relation which defines $A$ we
obtain $\ff^d g + \hh^d P(g) = 0$. Now $\gcd(g, P(g)) = 1$ in $K[S]$
since $P(0) \ne 0$, and so $g | \hh^d$. Since $\gcd(\ff, \hh )= 1$
we also have $\hh^d | g$. Write $g = \lambda \hh^d$ for some
$\lambda \in K^*$. We then have $\lambda \ff^d + P(\hh^d) = 0$, or
\begin{equation} \lambda \ff^d + \hh^{d}Q(\hh^d)=-P(0)
\label{EQ1}
\tag{$\dagger$}
\end{equation}
Since $P(0) \ne 0$ and $K$ is algebraically closed, (\ref{EQ1}) is
equivalent to a parametrization of the form $F^d + H^d Q(H) = 1$.
By
applying
Lemma~\ref{Lefschetz}
to
the hypothesis on
nonconstant parametrizations of this form, we must have
%
$\ff, \hh \in K^*$. So $g \in K^*$ as well which means $D(y) = 0$.
We now have $f = \ff \hh^{-1} h$ where $\ff \hh^{-1} \in K^*$. Now
$f(S) = x$ and $h(S) = z$ are elements of $A$, so $\ff \hh^{-1}$
belongs to the quotient field of $A$. By the way in which $K$ was
constructed (see the paragraph preceding Lemma \ref{PT1}), $\ff
\hh^{-1}$ must belong to the quotient field of $A^D$. Therefore we
can write $a f = b h$ for some elements $a,b \in A^D$, i.e. $ax =
bz$. But this relation combined with the relation $x^d y + z^d P(y)
= 0$ implies that $x$ and $z$ are algebraically dependent over
$A^D$. Then $x, z \in A^D$ by Lemma \ref{BS1}. But then $D =0$. This
contradiction completes the proof.
\end{proof}

\begin{remark}
In Proposition \ref{EX2} we must assume that no parametrizations of
the form $F^d + H^dQ(H) = 1$ exist.  If such parametrizations do
exist then it is unclear how to decide if $A$ in that proposition is
rigid.  For example, consider the ring
$A:=\KK[X,Y,Z]/(X^3Y+Z^3Y+Z^4)=\kk[x,y,z]$.  For any choice of
$\alpha, \ff, \hh \in K[s]$ we have that
\[ f=\alpha \ff (\ff^3+\hh^3), \quad g=-\alpha \hh^4, \quad  h=\alpha \hh(\ff^3+\hh^3) \]
is a solution to $f^3g+gh^3+h^4=0$ (perhaps assuming
$\gcd(\ff,\hh)=1$ if one likes the parametrization to be unique). It
is not clear from these parametrizations whether $A$ admits a
nontrivial LND.  However, the existence of
these parametrizations blocks us from using the techniques of this
section to decide rigidity. \end{remark}



\subsection{Rigidity of  $X^aY^b+Z^c+T^d$ for large $a,b,c,d$}

Write $A:=\C[X,Y,Z,T]/(X^aY^b+Z^c+T^d)$.

\begin{theorem}\label{abcd} Assume $a,b,c,d\geq 2$ and $a\geq b$. If $b^{-1}+c^{-1}+d^{-1}\leq 1$ then $A$ is rigid.
\end{theorem}

The proof is a direct application of
\begin{proposition} \label{doublemason}
Let $K$ be an algebraically closed field of characteristic 0.
Suppose $f,g,h,i\in K[S]$ and $f^ag^b+h^c+i^d=0$ for some positive
integers $a,b,c,d
$.
Assume $a\geq b$. If $b^{-1}+c^{-1}+d^{-1}\leq 1$ then $f,g,h,i$ are
all constant.
\end{proposition}

\begin{proof}
Write $F,G,H,I$ for the degrees of $f,g,h,i$, respectively. Using
Mason's Theorem \ref{Mason} for $n=3$ we have
$\max(aF+bG,cH, dI)<F+G+H+I$.
Suppose $\max(aF+bG,cH, dI) =aF+bG$. Now $aF+bG< F+G+H+I\leq
F+G+\frac{a}{c}F+\frac{b}{c}G+\frac{a}{d}F+\frac{b}{d}G$ and, since
$ a^{-1}+c^{-1}+d^{-1}\leq b^{-1}+c^{-1}+d^{-1}\leq 1$, this cannot
hold. We are left to assume $\max(aF+bG,cH, dI) = cH$ (without loss
of generality). Now $cH< F+G+H+I\leq
\frac{c}{a}H-\frac{b}{a}G+G+H+\frac{c}{d}H$ which means $H<
(a^{-1}+c^{-1}+d^{-1})H+(1-\frac{b}{a})G\leq
(a^{-1}+c^{-1}+d^{-1})H$. This too is not possible.
\end{proof}

\subsection{Rigidity of  $X_1^{d_1}+X_2^{d_2}+\ldots+X_n^{d_n}$ for large $d_i$}


\begin{lemma} \label{EX1} Let $n\geq 3$. Let $d_1,\ldots,d_n \in \N$ with all $d_i\geq 2$,
 and $\gcd(d_1,\ldots,d_n)=1$. If
$(d_1)^{-1}+(d_2)^{-1}+\ldots+(d_n)^{-1}\leq (n-2)^{-1}$ then
$\C[X_1,\ldots,X_n]/(X_1^{d_1}+X_2^{d_2}+\ldots+X_n^{d_n})$ is
rigid.
\end{lemma}

This is a direct application of Mason's Theorem \ref{Mason} and can
be found in \cite[Example 2.6]{FM07}.

\section{Filtrations}

%
%

Throughout this section $A$ will be a finitely generated
$\kk$-algebra.

\begin{definition} \label{D1}
A \textbf{semi-degree function} on $A$ is a function $\deg: A \to \Z
\cup \{-\infty\}$ satisfying
\begin{enumerate}
\item $\deg(a+b)\leq \max(\deg(a),\deg(b))$ for all $a, b \in A$,
\item $\deg(ab)\leq \deg(a)+ \deg(b)$ for all $a, b \in A$, and
\item $\deg(a)=-\infty$ if and only if $a = 0$.
\end{enumerate}
If equality always holds in condition 2 then $\deg$ is a
\textbf{degree function}. If  $\deg(c)=0$ for all $c\in \kk^*$, we
say that $\deg$ is a \textbf{$\kk$-(semi-)degree function}.
Since in this paper we are interested in $\kk$-algebras, we will
only consider $\kk$-(semi-)degree functions.
\end{definition}

\begin{definition}
A \textbf{filtration} $\F = \{\F_n(A)\}_{n \in \Z}$ of $A$ is a
chain
\[\cdots \subseteq\F_{n-1}(A)\subseteq \F_n(A) \subseteq \F_{n+1}(A)\subseteq \cdots \subseteq A \]
 of $\kk$-submodules of $A$ satisfying
\begin{enumerate}
\item $\F_{n}(A)\F_m(A)\subseteq \F_{n+m}(A)$ for all $n,m\in \Z$,
\item $\bigcup_{n \in \Z} \F_n(A) = A$, and
\item $\bigcap_{n \in \Z} \F_n(A) = 0$.
\end{enumerate}
\end{definition}

\begin{example}
If $A=\oplus_{i\in\N} A_i$  is a graded ring,
then $\F_n(A)=\oplus_{i\leq n} A_i$ defines a filtration of $A$.
\end{example}

\begin{remark}
A
$\kk$-semi-degree function $\deg$ on $A$ determines a filtration of
$A$ given by $\F_n(A):=\{a\in A~|~ \deg(a)\leq n \}$. Conversely, if
$\F = \{\F_n(A)\}$ is a filtration of $A$
 for which $c
\in \F_0(A) \backslash \F_{-1}(A)$ for all $c \in \kk^*$, then $\F$
determines a
$\kk$-semi-degree function on $A$ given by
$\deg(a) = \min\{n \in \Z ~|~ a \in \F_n(A)\}$. A
$\kk$-degree function on $A$ corresponds to such a filtration of $A$
with the additional property
\begin{enumerate}
\item[4.] If $a \in \F_n(A) \backslash \F_{n-1}(A)$ and $b \in \F_m(A)
\backslash \F_{m-1}(A)$, then $ab \in \F_{n+m}(A) \backslash
\F_{n+m-1}(A)$.
\end{enumerate}
\end{remark}



\begin{definition}
Let $\F$ be a filtration on $A$ with associated $\kk$-semi-degree
function $\deg$. Define, as a quotient of modules,
$\GR_{\F,n}(A):=\F_n(A)/\F_{n-1}(A)$ and $\GR_{\F}( A
):=\oplus_{n\in \Z}\GR_{\F,n}(A)$. Define the map $\gr:A\lp
\GR_{\F}(A)$ by $\gr(0) = 0$ and $\gr(a) =
a+\F_{\deg(a)-1}(A)$
if $a\not = 0$. $\GR_{\F}(A)$ is called the \textbf{graded ring
obtained from $A$ by $\F$} (or by $\deg$). Nonzero elements of
$\GR_{\F,n}(A)$ in $\GR(A)$ are called \textbf{homogeneous of degree
$n$}. When the filtration $\F$ is understood, we use the notation
$\GR_n(A)$ and $\GR(A)$ for simplicity.
\end{definition}

Note that every homogeneous element of degree $n$ in $\GR(A)$ is of
the form $\gr(a)$ for some $a \in \F_n(A) \backslash \F_{n-1}(A)$.
The set $\GR(A)$ is indeed a graded $\kk$-algebra. The linear
structure is clear. Multiplication of homogeneous elements is
defined as follows. If $\gr(a) = a+\F_{n-1}(A) \in \GR_{\F,n}(A)$
and $\gr(b) = b+\F_{m-1}(A)\in \GR_{\F,m}(A)$ are both nonzero, then
define $\gr(a) \gr(b)= ab+\F_{n+m-1}(A)$. This multiplication
can be extended to all elements of $\GR(A)$.
The following properties can be easily checked:
\begin{itemize}
\item $\gr(a)=0$ if and only if $a=0$.
\item $\gr(a+h)=\gr(a)$ if $\deg(a)>\deg(h)$.
\item $\gr(a+b)=\gr(a)+\gr(b)$ if $\deg(a)=\deg(b)=\deg(a+b)$.
\item $\gr(a^n)=\gr(a)^n$ if $\deg(a^n)=n\deg(a)$.
\item $\gr(ab)=\gr(a)\gr(b)$ if $\deg(ab)=\deg(a) + \deg(b)$.
\item $\GR(A)$ is reduced if and only if $\deg(a^n)=n \deg(a)$ for all $n\in
\N$.
\item $\GR(A)$ is a domain if and only if
$\deg(ab)=\deg(a) + \deg(b)$ for all $a,b\in A$, i.e. if $\deg$ is a
$\kk$-degree function on $A$.
\end{itemize}

When investigating the rigidity of a domain $A$, it is desirable for
$\GR(A)$ to be a domain.
But if $\deg(a)+ \deg(b)>\deg(ab)$ for some $a, b \in A$, it will
not be the case.

\begin{example}\label{Ex3}
Let $A:=\kk[X,Y]/(X^2-p(Y))$
, where $p\in \kk[Y]$. Write
$x:=X+(X^2-p(Y))$, so that $A=\kk[Y]\oplus x\kk[Y]$. Defining
$\deg(x)=1$ and $\deg(Y)=0$
indeed determines a $\kk$-degree function on $A$.
Now $\F_{n}(A)=\{0\}$ for $n < 0$, and $\F_{n}(A)=A$ for $n > 0$, while
$\F_{0}(A)=\kk[Y]$. Thus $\GR_1(A)=A/\kk[Y]$, $\GR_0(A)=\kk[Y]$, and
$\GR_{n}(A)=\{0\}$ for all other $n\in \Z$. Now
$\gr(x)=x+\F_{1-1}(A)=x+\kk[Y]$ is an element of $\GR_1(A)$. And
$\gr(x)^2=x^2+\F_{2-1}(A)=x^2+A= 0+A$, so $\gr(x)^2=0$. Note that
$\GR_1(A)=(\kk[Y]+x\kk[Y])/\kk[Y]=\gr(x)\kk[Y]$, so
$\GR(A)=\GR_0(A)\oplus \GR_1(A)= \kk[Y]\oplus \gr(x)\kk[Y]\cong
\kk[X,Y]/(X^2)$. Note that $\GR(A)$ is  the quotient of $\C[X,Y]$
modulo the ``highest degree part'' of $X^2-p(Y)$.

This
gives a funny example if $p(Y)=Y$. Then $A=\kk[X]$,
and $\GR(A)=\kk[X,Y]/(X^2)$.  Here $\kk[Y]$ is the image of the
original subring $\C[X^2]$, while the image of $X$ becomes
nilpotent.
\end{example}


Now we get to a pitfall.
It is tempting to believe that if $A$ is generated by
$x_1,\ldots,x_n$
then $\GR(A)$ will be generated by $\gr(x_1),\ldots,\gr(x_n)$. (In several papers this is not explicitly proven, but
true in that case.) In general this is not the case.

\begin{example}\label{Ex4}
Let $A=\kk[X,Y]$ and define $Z=Y+X^n$ for some $n > 1$. Define a
grading on $A$ by $\deg(a)=\deg_X(a)$ for all $a\in A$. Now $A$ is a
graded ring, and it is not hard to see that $\GR(A)$ is isomorphic
to $A$  and is generated by $\gr(X),\gr(Y)$. Now obviously $A$ is
generated by $X,Z$. But
$\gr(Z)=Y+X^n+\F_{n-1}(A)=X^n+\F_{n-1}(A)=\gr(X^n)=\gr(X)^n$, and
thus $\GR(A)$ is not generated by $\gr(X), \gr(Z)$.
\end{example}

While things can get quite messy, in many practical cases we can fortunately use the following
fact.


\begin{lemma} \label{l1}Let $A$
be a graded ring
generated by homogeneous elements $x_1, \ldots, x_n$.
Let $\F$ be the filtration defined by the grading of $A$. Then $
A\cong \GR_{\F}(A)$ by the isomorphism $x_i \mapsto
\gr(x_i)$, and thus $\GR_{\F}(A)$ is generated by
$\gr(x_1),\ldots,\gr(x_n)$.
\end{lemma}


\begin{proof}
We make use of the polynomial ring $\kk[X_1,\ldots,X_n]$, where the
$X_i$ are independent variables. Since $A$ is graded, we have for
any monomial $M(X_1,\ldots,X_n) \in \kk[X_1, \ldots, X_n]$ that
$M(\gr(x_1),\ldots,\gr(x_n))=\gr(M(x_1,\ldots,x_n))$. Define a
grading on
$\kk[X_1,\ldots,X_n]$
by making $X_i$ homogeneous of degree $\deg(x_i)$. Let $a\in A$ be
of degree $d$. Now $a=f(x_1,\ldots,x_n)$ for some $f\in
\kk[X_1,\ldots,X_n]$. Write $f=\sum f_i$ where $f_i$ is homogeneous
of degree $i\in \Z$. It is very well possible that there is some
$f_i\not = 0 $ for which $i> \deg(a)$. But in that case,
$f_i(x_1,\ldots,x_n)=0$ because the $x_i$ are homogeneous. So, we
could have taken $f-f_i$ instead of $f$. Thus we can assume that $f$
has degree $d$. Now
$\gr(a)=a+\F_{d-1}(A)=f(x_1,\ldots,x_n)+\F_{d-1}(A)
=f_d(x_1,\ldots,x_n)+\F_{d-1}(A)$ and since each monomial in $f_d$
is of degree $d$, we get that this equals
$f_d(\gr(x_1),\ldots,\gr(x_n))+\F_{d-1}(A)$.
\end{proof}

To emphasize the importance of starting with a graded ring:
If $A$ is not graded, then in the above argument it is not possible
to assume that the degree of $f$ is equal to
$\deg(a)$, and the equality
$f(x_1,\ldots,x_n)+\F_{d-1}(A)=f(\gr(x_1),\ldots,\gr(x_n))$ does not
hold. For instance, in Example \ref{Ex4} we have $\gr(Z-X^n)=Y + X^n
- X^n +\F_{-1}(A)=Y+\F_{-1}= \gr(Y)$, but
$\gr(Z)-\gr(X^n)=\gr(X^n)-\gr(X^n)=0$.

Lemma \ref{l1} is often applicable by the following trick. When
studying a ring $B$, view $B$ as a subring of a graded ring $A$ (for
which a particularly useful grading has been chosen). Using the
following lemma,
one can then attempt to find generators of $\GR(B)$.

\begin{lemma}
Let $A$ and $B$ be rings and let $\F$ be a filtration of $A$. If
$B\subseteq A$ then $\GR(B)\subseteq \GR(A)$.
\end{lemma}

\begin{proof} If $b \in B$ with $\deg(b) = d$, define
$\varphi(b+\F_{d-1}(B))=b+\F_{d-1}(A)$. This defines an injective
mapping $\GR_d(B) \lp \GR_d(A)$ which
extends to an injective ring homomorphism $\GR(B)\lp \GR(A)$.
\end{proof}

\section{Filtration techniques}

\begin{definition}
Let $A$ be a $\kk$-algebra with filtration $\F$. A linear map $L:A
\lp A$ \textbf{respects the filtration} if
\[
d_L:= \max_{a\in A, a\not = 0}\{ \deg(L(a))-\deg(a)\}
\]
is an integer. In this case we define a linear map $\gr(L): \GR(A)
\lp \GR(A)$ by $a + \F_{n-1}(A) \mapsto L(a) + \F_{n+d_L -1}(A)$ for
all $a \in \F_n(A) \backslash \F_{n-1}(A)$. By the choice of $d_L$,
note that if $L \ne 0$ then $\gr(L) \ne 0$.
\end{definition}

If $\deg(L(a))-\deg(a)=d_L$ for every $L(a)\not =0$, then $L$ can be called {\bf homogeneous of degree $d_L$}, and write $\deg(L):=d_L$. If $a$ has degree $e$, then $L(a)$ has degree $e+\deg(L)$ (or $-\infty$). In this paper, $L$ will always be a derivation.

\begin{lemma} \label{Disnonzero}
Let $A$ be a $\kk$-algebra with filtration $\F$. Let $D \in \LND(A)$
with $D \ne 0$. If $D$ respects the filtration on $A$ then $\gr(D)
\in \LND(\GR(A))$, with $\gr(D) \ne 0$,
%
and
$\gr(D)(\gr(a))=0$
for every $a \in A^D$.
Moreover, if $A$ is
finitely generated then $D$ respects the filtration.
\end{lemma}

\begin{proof}
It is straight forward to check the product rule and local
nilpotency of $\gr(D)$ for homogeneous elements on $\GR(A)$. These
properties then extend to all elements of $\GR(A)$, so that
$\gr(D)\in \LND(\GR(A))$.
It remains to show that $D$ respects the
filtration when $A$ is finitely generated. This follows from the
observation that if $A$ is generated by $a_1,\ldots, a_n$ then
\[
\max_{a\in A, a\not =0}\{ \deg(D(a))-\deg(a)\} = \max_{1\leq i\leq
p}\{ \deg(D(a_i))-\deg(a_i)\}.
\]
\end{proof}

\begin{corollary}
Let $A$ be a $\kk$-algebra with filtration $\F$. If $\GR(A)$ is rigid
then $A$ is rigid. 
%
\end{corollary}

In case $A:=\C^{[n]}/(f)$ and one has a grading on $\C^{[n]}$, then one can define a degree function on $A$ as follows:
for $g\in \C^{[n]}$, define $\deg(\bar{g})=\min\{ \grad(g+h)~|~h\in fA\}$. 
In this case, it is possible to describe $\GR(A)$ by the following result, which is a special case of
Proposition 4.1 of  \cite{KML07}.

\begin{proposition}\label{Prop:GrA}
Let $A:=\C^{[n]}/(f)$. Then $\GR_{\deg}(A)\cong \C^{[n]}/(\widehat{f}\,\,)$ where 
$\widehat{f}$ is the homogeneous highest degree part of $f$ with respect to the grading on $\C^{[n]}$.
\end{proposition}


\section{Examples using filtration techniques}

In this section we again consider examples of the form
$A:=\KK[X_1,\ldots,X_n]/(f)$.

%
%


\begin{lemma} \label{irreducible}
Suppose $A$ is a noetherian domain.
Let $a$ be an irreducible element of $A$.
If $A/aA$ is rigid, then the only  $D\in \LND(A)$ satisfying $D(a)=0$ is the zero derivation.
\end{lemma}



\begin{proof}
Let $D\in \LND(A)$ with $D(a)=0$. 
Suppose $D\not =0$.
Since $A$ is noetherian
there exists 
%
$D_1 \in \LND(A)$ and $b\in A^D$
such that $D=bD_1$
and $D_1(A)\not \subseteq aA$.
(See Lemma 2.2.14
of \cite{Thesis}.)
%
Now
$D_1(a)=0$, hence $D_1\mod{aA}$ is a well-defined
nonzero
LND on $A/aA$.
This gives a contradiction since $A/aA$ is rigid.
\end{proof}

As a remark, if $A$ is not noetherian one cannot conclude such things,
the main reason being that
a descending chain of principal ideals in a noetherian ring has an intersection which is principal.
 This is not true for non-noetherian rings as
Remark 2.2.16 of \cite{Thesis} shows.

\subsection{When is $X^a+Y^b+Z^c$ rigid?}

Remark that $X^a+Y^b+Z^c$ is always an irreducible polynomial for positive $a, b, c$. The following result is Lemma 4 in \cite{KZ00}.

\begin{theorem} \label{KalZai} The domain $A := \C[X,Y,Z]/(X^a+Y^b+Z^c)$ is rigid if both $a \geq 2$ and $b, c \geq 3$. 
\end{theorem}

Here is a very rough proof sketch. If $\frac{1}{a}+\frac{1}{b}+\frac{1}{c}\leq 1$ then use Mason's Theorem \ref{Mason}, and in the other cases use Makar-Limanov techniques (similar to in proofs in this article).

\subsection{When is $X^aY^b+Z^c+T^d$ a rigid domain? }

Remark that $X^aY^b+Z^c+T^d$ is irreducible for all positive $a,b,c,d$.
Consider the domain $A:=\C[X,Y,Z,T]/(X^aY^b+Z^c+T^d) =\C[x,y,z,t]$. First let us remove a few non-rigid cases. When $a=1$ then $cZ^{c-1}\partial_X -Y^b\partial_Z$ defines an LND on $A$. When $c=1$ 
then $A \cong \C^{[3]}$.
Also, if $c=d=2$ then
$A \cong \C[X,Y,Z,T]/(X^aY^b-ZT)$ which is not rigid.
The case $b = c = 2$ with $a$ even is also not rigid, as it admits an LND given by $dT^{d-1} \partial_Y - d i X^{a/2}T^{d-1} \partial_Z - 2X^{a/2}(X^{a/2} Y - i Z) \partial_T$, where $i^2 = -1$. (Remark that $X$ and $X^{a/2} Y - iZ$ are in the kernel of this derivation. This example was shown to us by Gene Freudenburg.)
Hence we assume $a,b,c,d\geq 2$ and not $c=d=2$, and we assume not $b=c=2$ if $a$ is even.
We suspect that these assumptions ensure $A$ to be rigid, but a couple of cases elude us.
In the next few lemmas we refer to the following technical assumption:

\begin{assumption} \label{assumption} Both
(1) no more than one of $a,c,d$ equals 2, and
(2) no more than one of $b,c,d$ equals 2.
\end{assumption}


\begin{lemma} \label{3fold.L2}

Suppose $D\in \LND(A)$. If \ref{assumption} holds, and if one of $x,y,z,t$ is in $A^D$, then $D=0$. 
\end{lemma}

\begin{proof}
If $x$ (or $y$) $\in A^D$, then 
$D$ induces an LND on
%
$\C[Y,Z,T]/(\lambda Y^b+Z^c+T^d)\cong \C[Y,Z,T]/(Y^b+Z^c+T^d)$
(where $ \lambda \in \C^*$ is the image of $x^a \in A$).
%
This ring is rigid by Theorem \ref{KalZai},
%
so $D=0$
by Lemma \ref{irreducible}. 
If $z$ (or $t$) $\in A^D$, then
$D$ induces an LND on
%
$B:= A/(z) \cong C[X,Y,T]/(X^a Y^b+ T^d)$.
%
By Theorem \ref{case1} this ring is rigid, 
%
so
$D=0$
by Lemma \ref{irreducible}. 
\end{proof}

Define two gradings on $A$, according to the weight distributions
$\deg_1(x,y,z,t)=(cd,0,ad,ac)$ and  $\deg_2(x,y,z,t)=(0,cd,bd,bc)$. Assuming $A$ is not rigid, using Lemma \ref{Disnonzero} we
have a nonzero $D\in \LND(A)$ which is homogeneous with respect to both gradings.
We will use this homogeneous LND in the next few lemmas. By Lemma \ref{Disnonzero} we can find a homogeneous element in $A^D$.

\begin{lemma} \label{3fold.L3}
If \ref{assumption} holds and if $\gcd(a,b)=1$ then $A$ is rigid.
\end{lemma}

\begin{proof}
%
Suppose $\deg_i (x^p y^q z^r t^s) = \deg_i(x^P y^Q z^R t^S)$, $i=1,2$, for some monomials. Then $cd(p-P) + ad(r-R) + ac(s-S) = 0$ and $cd(q-Q) + bd(r-R) + bc(s-S)=0$. From this we find $b(p-P) = a(q-Q)$. Since $\gcd(a,b)=1$, we conclude $p-P = al$ and $q-Q = bl$ for some $l \in \Z$. So $x^p y^q z^r t^s = x^P y^Q (x^a y^b)^v z^r t^s = x^P y^Q (-z^c - t^d)^v z^r t^s$. Therefore, homogeneous elements can be written as $x^i y^j \tilde{f}(z,t)$.

Suppose now that $f = x^i y^j \tilde{f}(z,t)$ is a homogeneous element in $A^D$. If $i>0$ or $j>0$, then  $D(x)=0$ or $D(y)=0$ 
by Lemma \ref{BS1} part b)
and we are done by Lemma \ref{3fold.L2}. 
So we assume that
%
$f \in
\C[z,t]$. Since $\trdeg_{\C} A^D\geq
2$ by Lemma \ref{BS1} part e), we find two algebraically independent elements $f,g \in \C[z,t]$ such that $D(f)=D(g)=0$.
But since $z$ is algebraically dependent of these $f,g$ we also have $D(z)=0$ by Lemma \ref{BS1} part c). By Lemma \ref{3fold.L2} we are done.
\end{proof}

\begin{lemma} \label{3fold.L4}
If \ref{assumption} holds and if  $\gcd(ac,d)=1$ then $A$ is rigid. 
\end{lemma}

\begin{proof}
Let $f=\sum_{i,j,k,l} f_{i,j,k,l} x^iy^jz^kt^l \in A^D$ be homogeneous w.r.t.
$\deg_1$. We may assume that $f_{i,j,k,l}=0$ if $l\geq d$. 
Using $\gcd(a,d)=\gcd(c,d)=1$, a computation (similar to that in Lemma \ref{3fold.L3}) yields that $\deg_1(x^iy^jz^kt^l)=\deg_1(x^Iy^Jz^Kt^L)$ only if $l-L$ is a multiple of $d$. This tells us that $f=t^m \tilde{f}(x,y,z)$ for some $m$.
If $m>0$ then $D(t)=0$ and Lemma \ref{3fold.L2}
shows that $D=0$. Thus we may assume 
$f \in
\C[x,y,z]$.
Since $\trdeg_{\C} A^D\geq 2$, we can find two algebraically independent homogeneous elements $f, g \in \C[x,y,z] \cap A^D$.

Another computation on degrees shows that
%
%
if
$h\in\C[x,y,z]$
is homogeneous w.r.t both degrees, then $h$
can be written as
\[ h=x^py^qz^r\sum_{i=0}^n h_i (x^{\alpha}y^{\beta})^i(z^{\gamma})^{n-i} =  x^py^qz^r\tilde{h}\prod_{i=1}^n(x^{\alpha}y^{\beta}+\tilde{h}_iz^{\gamma}) \]
where $h_i, \tilde{h}, \tilde{h}_i\in \C$, $p,q,r\in \N$, and $\alpha=a/e, \beta=b/e, \gamma = c/e$ where $e=\gcd(a,b,c)$.
Since $f,g$ are algebraically independent and homogeneous, this means that there are at least two algebraically independent elements of the form
$x,y,z$, or $x^{\alpha}y^{\beta}+\lambda z^{\gamma}$ in $A^D$. Using Lemma \ref{3fold.L2} we may assume that $x,y,z\not \in A^D$. Hence we have two kernel elements of the form 
$x^{\alpha}y^{\beta}+\lambda z^{\gamma}$ and $x^{\alpha}y^{\beta}+\mu z^{\gamma}$ where $\lambda\not = \mu$.  But then their difference
is also in $A^D$, from which follows that $z\in A^D$. So $D=0$.
\end{proof}

\begin{lemma} \label{3fold.L5}
If \ref{assumption} holds and if  $\gcd(a,cd)=\gcd(b,cd)=1$ then $A$ is rigid.
\end{lemma}

\begin{proof}
%
Suppose
$f \in A^D$
is homogeneous with respect to both degrees.
A computation yields that $f=x^iy^j\tilde{f}(z,t)$. We may hence assume that we have two  algebraically independent $f, g\in \C[z,t]
\cap A^D$.
Since $z$ is algebraically dependent of $f$ and $g$, $D(z)=0$ and hence we are done.
\end{proof}

\begin{lemma} \label{3fold.L6}
If $a \geq 3$ is odd and $d \geq 3$ then $X^a Y^2+Z^2 + T^{d}$ is rigid.
\end{lemma}


\begin{proof}
As shown in Lemma \ref{3fold.L3}, homogeneous elements of $A$ can be written as $x^i y^j f(z,t)$. Furthermore, $f(z,t)$ is a product of factors each of the form $ \mu z^2+ \lambda t^d$ (if $d$ is odd) or $ \mu z+\lambda t^{d/2}$ if $d$ is even.
%
%
If one of $y$, $z$, or $t$ is in $A^D$ then arguing as in Lemma \ref{3fold.L2} we find that $D=0$.

Since there are at least two algebraically independent elements in $A^D$, we can conclude that either
(1) two algebraically independent homogeneous elements $f(z,t)$, $g(z,t) \in A^D$, or (2) $x$ and one homogeneous element $f(z,t)$ are in $A^D$.
Case (1) implies, since $z$ and $t$ are  algebraically dependent of $f,g$, that $D=0$. 
If $d$ is odd then case (2) is solved as $\mu z^2+\lambda t^d \in A^D$ implies by Lemma \ref{MiniMason} that $t\in A^D$ or $z\in A^D$ and so $D=0$.

The only case left is if $d=2e$. We have $x, z - \lambda t^e\in A^D$ for some $\lambda \in \C$. 
Then $\bar{D}:=D \mod (z-\lambda t^e)$ is an LND on  $X^a Y^2+(\lambda^2+1) T^{d}$. If $\lambda^2+1\not =0$ then $\bar{D}=0$ by Theorem \ref{case1}, and 
$D=0$
by Lemma \ref{irreducible}.
Thus we are left to assume that $\lambda^2 = -1$.
Now 
define $\bar{D}:=D \mod (z- \lambda t^e-1)$, which is an LND on $X^a Y^2+(1+2\lambda T^e)$. This surface is rigid, as its top degree part $X^{a}+ T^e$ is rigid 
by Theorem \ref{case1}. Thus $\bar{D}=0$, and hence $D=0$
again by Lemma \ref{irreducible}.
\end{proof}


\begin{theorem} \label{abcdTHM}
With the exceptions of the  cases mentioned in Remark \ref{Leftover} below, $X^aY^b+Z^c+T^d$ is rigid. 
\end{theorem}

\begin{proof}
We assume $a\geq b$ and  $d\geq c$.
We can apply Theorem \ref{abcd} to remove many cases. We are left with only a few possibilities for $\{b,c,d\}$, with $\gcd(a,b) > 1$.\\
Case $\{2,3,5\}$ is completely covered by Lemma \ref{3fold.L4}.\\
Case $\{2,3,4\}$. Lemma \ref{3fold.L4} works if $c=3$ or $d=3$
(interchanging roles of $a$ and $b$ in that Lemma as necessary). If $b=3$ then we can apply Lemma \ref{3fold.L5} unless $\gcd(a,4,2)>1$ (i.e. $a$ is even), and we can apply Lemma \ref{3fold.L3} unless $\gcd(a,3)>1$. Thus left is
Remark \ref{Leftover} (i).\\
Case $\{2,3,3\}$. Lemma \ref{3fold.L4} works if $c=2$. If $b=2$ then we can use Lemma \ref{3fold.L3} unless $a$ is even, and use Lemma \ref{3fold.L5} unless $a=3\lambda$, hence left is Remark \ref{Leftover} (ii).\\
Case $\{2,2,m\}$, where $m\geq 3$. We exclude $b=m$ (for not being rigid)
and so we assume $b=c=2$ and $d \geq 3$. We assume $a \geq 3$ is odd, since otherwise $A$ is not rigid.
The remaining case is covered by Lemma \ref{3fold.L6}.
%
\end{proof}

\begin{remark} \label{remark}\label{Leftover}
If $a,b,c$ or $d$ equals 1, or if $c=d=2$, 
or if $b = c = 2$ with $a$ even,
then $X^aY^b+Z^c+T^d$ is not rigid. The following cases are unproven to be rigid.
\[ \begin{array}{l r r}
\textup{(i)~}& X^{6a}Y^3+Z^2+T^4, & \quad a\geq 1\\
\textup{(ii)~}& X^{6a}Y^2+Z^3+T^3, & \quad a\geq 1
\end{array}\]

\end{remark}

\subsection{When is $X^a+Y^b+Z^c+T^d$ rigid?}

Write $A:=\C[x,y,z,t]=\C[X,Y,Z,T]/(X^a+Y^b+Z^c+T^d)$,  where we  assume $a,b,c,d\geq 2$ and at most one equals $2$.
(If one of $a,b,c,d$ equals 1 then $A \cong \C^{[3]}$, where if for example $a=b=2$ we have a non-rigid ring
isomorphic to $XY+Z^c+T^d$.
While Theorem \ref{KalZai} classifies the surface $X^a+Y^b+Z^c = 0$,
it turns out that a complete classification for the one-dimension higher case
 $X^a+Y^b+Z^c+T^d=0$ is out of reach for now. 
It will turn out that
many
cases will not be covered by the results from this section or Lemma \ref{EX1}.

Define a
degree function on $A$ by $\deg(x,y,z,t)=(bcd,acd,abd,abc)$.

\begin{lemma} \label{CB2}  
Let $D\in \LND(A)$. If one of $x,y,z,t\in A^D$, then $D=0$.
\end{lemma}

\begin{proof} 
Suppose without loss of generality that $D(t)=0$. Then
$A/(t)$ is a domain 
and is rigid by Theorem \ref{KalZai}. So $D \mod tA =0$, hence using Lemma \ref{irreducible}, $D=0$.
\end{proof}

\begin{lemma}\label{CB1} 
Suppose that $\gcd(abc,d)=1$. If $f\in A$ is homogeneous w.r.t. $\deg$, then
$f=t^ig$ for some $i\in \N$ and some $g\in \C[x,y,z]$.
\end{lemma}

\begin{proof}
This is seen by a calculation similar to those in Lemmas \ref{3fold.L3}, \ref{3fold.L4}, and \ref{3fold.L5}.
\end{proof}

\begin{lemma} \label{CB3}
Suppose that $\gcd(abc,d)=1$ and that $A$ is not rigid. Then there exists some nonzero $D\in \LND(A)$
of the form
\[ D=dt^{d-1}\delta - Q(x,y,z) \frac{\partial}{\partial t} \]
where $\delta\in \DER(\C[x,y,z])$, $Q=\delta(x^a+y^b+z^c)$, and $D(Q)=\delta(Q)=0$. 
\end{lemma}

\begin{proof}
We can assume $D$ is homogeneous, and thus $D(t)$ is homogeneous. By Lemma \ref{CB1} we know that $D(t)=t^i p(x,y,z)$ for some $i$. If $i>0$ then $D(t)=0$ by Lemma 
\ref{BS1} part b), and Lemma \ref{CB2} tells us that $D=0$. So $D(t)=-Q(x,y,z)$ for some $Q$. 
$D$ sends homogeneous elements to homogeneous elements, and there is a constant $e$ such that $D(A_p)\subseteq A_{p+e}$. Since $\deg(D(t))=\deg(Q)=0 \mod(d)$,
$e=\deg(D(t))-\deg(t)=-abc \mod(d)$. Thus, $\deg(D(x))= \deg(x)+e= -abc \mod(d)$. 
Thus $D(x)\in t^{d-1}\C[x,y,z]$ because of Lemma \ref{CB1}. The same holds for $D(y)$ and $D(z)$, hence $D=dt^{d-1}\delta -Q\partial_t$ for some $\delta\in \DER(\C[x,y,z])$. 
Now $dt^{d-1}\delta(x^a+y^b+z^c)=D(x^a+y^b+z^c)=D(-t^d)=dt^{d-1}Q$, hence
$\delta(x^a+y^b+z^c)=Q$. Furthermore, since $D^2(t)=D(Q)=dt^{d-1}\delta(Q)$
we get using lemma \ref{BS1} part d) that $\delta(Q)=0$ and hence $D(Q)=0$.
\end{proof}

Unfortunately, the above explicit description was not enough to prove rigidity. We mention it to aid future work.
%
With a couple more assumptions we can salvage the following result.



\begin{theorem} \label{CB4}
If $\gcd(ab,c)=\gcd(abc,d)=1$ and if $\gcd(a,b) \ne a$ (or $b$), then $A$ is rigid. 
\end{theorem}

\begin{proof}
We may assume  $D\in \LND(A)$ to be homogeneous. 
By our assumptions, homogeneous elements in $A$ now have the form $x^i y^j z^k t^l f(x^{\alpha}, y^{\beta}, z^c)$, where $\alpha = a/\gcd(a,b) > 1$ and $\beta = b/\gcd(a,b) > 1$.
Since $\trdeg(A^D)\geq \trdeg(A)-1=2$, there exist two algebraically independent homogeneous elements $f,g
\in A^D$. 
By Lemma \ref{CB2} we can assume $f$ and $g$ are polynomials in $x^{\alpha}$, $y^{\beta}$, $z^c$. By Lemma 11.2 in \cite{KML07} two of the three elements $x$, $y$, and $z$ are in $A^D$. So $D=0$.
\end{proof}

\subsection{A cautionary tale: $D(X^2+Y^3+Z^5)=0$ does not imply $D(X)=D(Y)=D(Z)=0$}

As a remark, we want to state that one should be careful when using 
techniques
discussed in this paper.
If
$A$ is a domain over $\C$ in which there are pairwise algebraically independent elements $x,y,z$ such that $x^2+y^3+z^5=0$, it does not necessarily mean that $x,y,z\in A^D$ for any $D\in \LND(A)$. Consider the following example.
Define
\[
R:=\C[f,g,h]:=\C[F,G,H]/(F^a+G^b+H^c)
\]
Let
%
$A \subset R[S]$ be defined by $A:= R[S^if, S^ig, S^ih ~|~ i\geq 1]$.
Define $x:=S^{bc}f$, $y:=S^{ac}g$, and $z:=S^{ab}h$.
Let $D:=\partial_S$.
Then $D(x^a+y^b+z^c)=0$ but not $D(x)=0$.
The point is that, though $x,y,z$ are pairwise algebraically independent, they are not generators of the ring $A$. (Note that the example can be modified to get $A$ to be a finitely generated $\kk$-algebra.)

\section{Open questions}

Let us mention a few explicit rings for which rigidity is unsettled.
The first is a surface in $\C^3$. The rest are 
hypersurfaces in $\C^4$.
%
Are
any of
them
rigid?
\[\begin{array}{l}
X^3Y+Z^3Y+Z^4\\
X^6Y^3+Z^2+T^4\\
X^6Y^2+Z^3+T^3\\
X^2+Y^3+Z^3+T^3\\
X^3+Y^3+Z^3+T^3\\
X^2+Y^3+Z^5+T^{15}\\
\end{array}\]

\section*{Acknowledgement}

The authors are grateful to Gene Freudenburg for providing several useful comments.

\end{document}